\documentclass[12pt]{amsart}
\usepackage[margin=1in]{geometry}
\usepackage{amsmath,amssymb,amsthm,mathtools}
\usepackage[T1]{fontenc}
\usepackage{lmodern, microtype}
\usepackage[hidelinks]{hyperref}
\newtheorem{theorem}{Theorem}[section]
\newtheorem{lemma}[theorem]{Lemma}
\newtheorem{proposition}[theorem]{Proposition}
\newtheorem{corollary}[theorem]{Corollary}
\theoremstyle{remark}
\newtheorem{remark}[theorem]{Remark}
\numberwithin{equation}{section}

\newcommand{\R}{\mathbb{R}}
\newcommand{\Z}{\mathbb{Z}}
\newcommand{\EE}{\mathbb{E}}
\newcommand{\HH}{\mathcal{H}}
\newcommand{\eps}{\varepsilon}
\newcommand{\dimH}{\dim_{\mathrm H}}
\newcommand{\dimB}{\overline{\dim}_{\mathrm B}}
\newcommand{\dimT}{\dim_{\mathrm{top}}}
\DeclareMathOperator{\Holdim}{H\ddot{o}ldim}
\DeclareMathOperator{\Cdim}{Cdim}
\DeclareMathOperator{\diam}{diam}

\title[Dimension reduction by bi-H\"older maps]
{Dimension reduction by bi-H\"older maps}

\author{Leonid V. Kovalev}
\address{215 Carnegie, Department of Mathematics, Syracuse University,
Syracuse, NY 13244, USA}
\email{lvkovale@syr.edu}

\subjclass[2020]{Primary 28A78; Secondary 28A80, 30L05, 54F45}

\keywords{H\"older dimension, Hausdorff dimension, box-counting dimension, Minkowski dimension, 
bi-H\"older embedding, singular function}

\begin{document}

\begin{abstract}
The H\"older dimension of a metric space is the infimum of its Hausdorff dimension over all bi-H\"older changes of the metric. Colvin asked whether it can take noninteger values. 
We show that every metric space of finite upper box dimension has integer H\"older dimension. 
For a bounded set $E\subset\R^n$ the dimension reduction is carried out by a homeomorphism of $\R^n$, built from randomly translated singular functions, which brings the Hausdorff dimension of $E$ arbitrarily close to its integer part.
\end{abstract}

\maketitle

\section{Introduction}\label{sec-intro}

For a nonempty bounded metric space $(X,d)$, its \emph{H\"older dimension}
is
\begin{equation}\label{eq-HD}
\Holdim X=\inf\{\dimH(X,\rho)\colon c\,d^\alpha\le\rho\le C\,d^\beta
\text{ for some } c,C,\alpha,\beta>0\}.
\end{equation}
where $\dimH$ is the Hausdorff dimension.
The metrics $\rho$ in~\eqref{eq-HD} will be called \emph{admissible}. They
are the pullbacks of metrics under bi-H\"older homeomorphisms of $X$ onto
other metric spaces. 
Topological dimension $\dimT$ is invariant under
homeomorphisms and does not exceed Hausdorff dimension, by a theorem of
Szpilrajn \cite[Theorem VII 2]{HurewiczWallman}. Hence
\begin{equation}\label{eq-basic}
\dimT X\le\Holdim X\le\dimH X.
\end{equation}

H\"older dimension was introduced by Colvin~\cite{Colvin}, who asked
whether it can take noninteger values \cite[Question 1.14]{Colvin}. We
answer this question for spaces of finite upper box dimension. 

For a metric space $X$ and $r>0$, let $N_X(r)$ be the least number of balls of radius $r$
centered in $X$ that cover $X$, with $N_X(r)=\infty$ if there is no finite cover. The \emph{upper box dimension} of $X$ is
\[
\dimB X=\limsup_{r\to0}\frac{\log N_X(r)}{\log(1/r)}.
\]
A space with $\dimB X<\infty$ is totally bounded. All bounded subsets of Euclidean spaces, and all bounded doubling metric spaces, have finite upper box dimension. 

\begin{theorem}\label{thm-integer}
Let $(X,d)$ be a metric space with $\dimB X<\infty$. Then $\Holdim X$ is an integer and therefore 
\begin{equation}\label{eq-integer}
\dimT X\le\Holdim X\le\lfloor\dimH X\rfloor.
\end{equation}
Moreover, the infimum in~\eqref{eq-HD} is unchanged if it is restricted to
bi-H\"older embeddings of $X$ into $\mathbb R^n$, provided that $n$ is allowed to vary.
\end{theorem}

The Euclidean ingredient of the proof is a dimension reduction that takes place
within the ambient space. We equip $\R^n$ with the maximum norm; by the
equivalence of norms, the statements below hold for the Euclidean norm as
well.

\begin{theorem}\label{thm-euclidean}
Let $k\ge0$ be an integer and let $E\subset\R^n$ be a bounded set with
$\dimH E<k+1$. For every $\eps>0$ there is a homeomorphism
$F(x)=(f_1(x_1),\dots,f_n(x_n))$ of $\R^n$, with increasing coordinate
functions and bi-H\"older on every bounded set, such that $\HH^{k+\eps}(F(E))=0$. 
In particular, $\dimH F(E)\le\lfloor\dimH E\rfloor+\eps$.
\end{theorem}

Colvin proved that the H\"older dimension of a compact locally
self-similar space is equal to its topological dimension, and that the
H\"older dimension of a compact doubling space does not exceed its capacity
dimension \cite[Corollaries 1.8 and 1.10]{Colvin}. He also constructed a Cantor set on the line with
H\"older dimension $1$ \cite[Theorem 8.1]{Colvin}, so the lower bound
in~\eqref{eq-integer} can be strict. The upper bound can be strict as well: e.g., the space 
$X = [0, 1]$ with the metric $d_X(x,y)=|x-y|^{1/2}$ has Hausdorff dimension $2$ and H\"older dimension $1$.

It is instructive to compare H\"older dimension with conformal dimension, defined as
the infimum of Hausdorff dimension over all quasisymmetric images~\cite{MackayTyson,Pansu}. 
Quasisymmetric maps between bounded uniformly
perfect spaces are bi-H\"older \cite[Corollary 11.5]{Heinonen}, so for such
spaces $\Holdim X\le\Cdim X$, as observed in~\cite{Colvin}. 
Moreover, the quasisymmetric dimension reduction of~\cite{Kovalev} implies that $\Holdim X \notin (0,1)$ for any bounded
metric space $X$~\cite[Theorem 3.3.5]{ColvinThesis}. 
To give an example with $\Holdim X < \Cdim X$, let $C_t\subset[0,1]$ be the self-similar Cantor set of Hausdorff dimension $t\in (0, 1)$. Then $\Cdim(C_t\times [0,1])=1+t$ by~\cite[Proposition 4.1.11]{MackayTyson} while $\Holdim(C_t\times[0,1])=1$ because all spaces 
$C_t\times [0, 1]$ are bi-H\"older equivalent.

The proof of Theorem~\ref{thm-euclidean} changes each coordinate by an
independent random translate of a singular increasing function. 
Theorem~\ref{thm-integer} follows by combining
Theorem~\ref{thm-euclidean} with the embedding theorem of Foias and Olson
for spaces of finite upper box dimension (Theorem~\ref{thm-embedding}).

Singular changes of variable that lower dimension go back to Tukia~\cite{Tukia}, who constructed a quasisymmetric self-map of $[0,1]$
mapping a set of full measure onto a set of arbitrarily small dimension. 
More recently, $n$-dimensional singular quasisymmetric maps were constructed by Romney~\cite{Romney}. 
These results control a set of full measure. Theorem~\ref{thm-euclidean} controls an arbitrary set, at the
cost of stopping at the integer part and of using bi-H\"older rather than quasisymmetric maps. It can also be read as a
statement about almost every translate of a set under one fixed map, in the spirit of~\cite{BMT,BHW}.

Whether every bounded metric space of finite Hausdorff dimension has
integer H\"older dimension remains open.

\section{A random singular change of variable}\label{sec-scalar}

Fix $0<p<1/2$ and let $\mu$ be the Bernoulli probability measure on
$[0,1]$ that assigns to the left and right halves of each dyadic interval
the fractions $p$ and $1-p$ of its measure. Its distribution function
$h$ is Salem's singular function~\cite{Salem}; since $\mu$ has no
atoms and has full support, $h$ is an increasing homeomorphism of $[0,1]$.
Extend it to $\R$ by $h(s+m)=h(s)+m$ for $m\in\Z$. For $0<q<1$ put
\begin{equation}\label{eq-parameters}
a=-\log_2 p,\qquad b=-\log_2(1-p),\qquad
\tau=\tau_q(p)=1-\log_2\bigl(p^q+(1-p)^q\bigr).
\end{equation}
Then $a>1>b>0$ and $q<\tau<1$. For fixed $q$ we have $\tau_q(p)\to1$ as
$p\to0$.

\begin{lemma}\label{lem-scalar}
There are constants $c_p,C_p>0$ such that
\begin{equation}\label{eq-power}
c_p|s-t|^a\le|h(s)-h(t)|\le C_p|s-t|^b\qquad(|s-t|\le1).
\end{equation}
If $T$ is uniformly distributed on $[0,1]$, then for every $u\in\R$ and
every $0<r\le1$,
\begin{equation}\label{eq-moments}
\EE\bigl[h(u+T+r)-h(u+T)\bigr]=r,\qquad
\EE\bigl[(h(u+T+r)-h(u+T))^q\bigr]\le4r^\tau.
\end{equation}
\end{lemma}

\begin{proof}
A dyadic interval of level $j$ has $\mu$-measure between $p^j$ and
$(1-p)^j$. An interval of length $r\le1$ contains a dyadic interval of
length at least $r/4$, and is contained in the union of at most two dyadic
intervals of length less than $2r$. This proves~\eqref{eq-power}.

The increment in~\eqref{eq-moments} is the measure of a translated interval
of length $r$ under the periodic extension of $\mu$, so its expectation
is $r$ by Fubini's theorem. For the second estimate choose a dyadic length
$l=2^{-j}$ with $r\le l<2r$. On the circle of length $1$, a translated
interval of length $r$ meets a fixed dyadic interval of length $l$ for a
set of translations of measure at most $r+l\le2l$. Since $q<1$,
subadditivity of $x\mapsto x^q$ and the binomial formula give
\[
\EE\bigl[(h(u+T+r)-h(u+T))^q\bigr]
\le 2l\sum_{|I|=2^{-j}}\mu(I)^q
=2l\bigl(p^q+(1-p)^q\bigr)^j=2l^\tau\le4r^\tau,
\]
the sum being over the dyadic intervals of level $j$. 
\end{proof}

For a random vector $T=(T_1,\dots,T_n)$ with independent coordinates
uniformly distributed on $[0,1]$, define
\begin{equation}\label{eq-random-map}
F_T(x)_i=h(x_i+T_i)-h(T_i),\qquad 1\le i\le n.
\end{equation}
Each $F_T$ is a homeomorphism of $\R^n$, and~\eqref{eq-power}
gives
\begin{equation}\label{eq-ambient-power}
c_p\|x-y\|_\infty^a\le\|F_T(x)-F_T(y)\|_\infty\le C_p\|x-y\|_\infty^b
\qquad(\|x-y\|_\infty\le1),
\end{equation}
with constants independent of $T$. Since $F_T(x+m)=F_T(x)+m$ for
$m\in\Z^n$, these estimates imply that $F_T$ is bi-H\"older on every
bounded set, with constants depending on its diameter.

\section{Rectangular covers and Euclidean compression}
\label{sec-compression}

For $s\ge0$ and $A\subset\R^n$, let
$\HH^s_\infty(A)=\inf\sum_i(\diam A_i)^s$ be the $s$-dimensional Hausdorff
content of $A$, the infimum being taken over all countable covers
$\{A_i\}$ of $A$. Recall that $\HH^s(A)=0$ if and only if
$\HH^s_\infty(A)=0$. In the maximum norm, a cube of side length $r$ has
diameter $r$.

Let $Q=\prod_{i=1}^n[u_i,u_i+r]$ be a cube with side length $0<r\le1$. Its
image under~\eqref{eq-random-map} is a rectangle with independent positive
side lengths $L_i=h(u_i+T_i+r)-h(u_i+T_i)$, which by
Lemma~\ref{lem-scalar} satisfy
\begin{equation}\label{eq-sides}
\EE L_i=r,\qquad \EE L_i^q\le4r^\tau.
\end{equation}

\begin{lemma}\label{lem-rectangle}
Let $0\le k<n$ and $0<q<1$. Then
\begin{equation}\label{eq-product-bound}
\HH^{k+q}_\infty(F_T(Q))\le\Psi_Q(T) := 2^k\sum_{i=1}^n\
\sum_{\substack{J\subset\{1,\dots,n\}\setminus\{i\}\\ |J|=k}}
L_i^q\prod_{j\in J}L_j ,
\end{equation}
and
\begin{equation}\label{eq-expected-cost}
\EE\,\Psi_Q(T)\le C_{n,k}\,r^{k+\tau}.
\end{equation}
\end{lemma}

\begin{proof}
Arrange the side lengths of $F_T(Q)$ in decreasing order,
$L_{(1)}\ge\dots\ge L_{(n)}$, and let $t=L_{(k+1)}$. Dividing each of the
$k$ longest sides into intervals of length at most $t$ covers $F_T(Q)$ by
at most
\[
\prod_{j=1}^k\left\lceil\frac{L_{(j)}}{t}\right\rceil
\le2^kt^{-k}\prod_{j=1}^kL_{(j)}
\]
cubes of side $t$; for $k=0$ this is a single cube. Hence
$\HH^{k+q}_\infty(F_T(Q))\le2^kL_{(k+1)}^q\prod_{j=1}^kL_{(j)}
\le\Psi_Q(T)$. Each product in~\eqref{eq-product-bound} involves distinct,
hence independent, side lengths, so by~\eqref{eq-sides} its expectation is
at most $4r^{k+\tau}$. There are $n\binom{n-1}{k}$ such products.
\end{proof}

\begin{proof}[Proof of Theorem~\ref{thm-euclidean}]
If $k\ge n$, there is nothing to prove. Assume $k<n$ and choose 
$0<q<\min\{1,\eps\}$. Since $\tau_q(p)\to1$ as $p\to0$, we can choose $p$ so that $\dimH E<k+\tau$, where $\tau=\tau_q(p)$. Then
$\HH^{k+\tau}_\infty(E)=0$. Enclosing the members of suitable covers in cubes of slightly larger positive side, we obtain for every $m\ge1$ a countable cover $\mathcal Q_m$ of $E$ by cubes with side lengths 
$0<r_Q<1$ such that $\sum_{Q\in\mathcal Q_m}r_Q^{k+\tau}<2^{-m}$. Let
\[
Z_m(T)=\sum_{Q\in\mathcal Q_m}\Psi_Q(T).
\]
By countable subadditivity of $\HH^{k+q}_\infty$ and~\eqref{eq-product-bound},
$\HH^{k+q}_\infty(F_T(E))\le Z_m(T)$ for every $m$. By monotone convergence
and~\eqref{eq-expected-cost}, $\EE Z_m\le C_{n,k}2^{-m}$. Hence
$Z=\inf_mZ_m$ satisfies $\EE Z=0$, and so $Z(T)=0$ for almost every $T$.
For such $T$ we have $\HH^{k+q}_\infty(F_T(E))=0$, hence
$\HH^{k+q}(F_T(E))=0$.
\end{proof}

Theorem~\ref{thm-euclidean} also gives integrality for two infima taken
within a fixed ambient space.

\begin{corollary}\label{cor-restricted}
Let $E\subset\R^n$ be a bounded set. Each of the infima
\begin{itemize}
\item[(i)] $\inf\{\dimH f(E)\colon f\colon E\to\R^n$ is a bi-H\"older
embedding$\}$,
\item[(ii)] $\inf\{\dimH F(E)\colon F\colon\R^n\to\R^n$ is a homeomorphism
that is bi-H\"older on bounded sets$\}$
\end{itemize}
is an integer between $\dimT E$ and $\lfloor\dimH E\rfloor$.
\end{corollary}

\begin{proof}
Both classes of maps contain the identity and are closed under
postcomposition with the maps of Theorem~\ref{thm-euclidean}. If one of the
infima were strictly between integers $k$ and $k+1$, some image $f(E)$ or
$F(E)$ of dimension less than $k+1$ could be compressed to dimension at
most $k+\eps$ for every $\eps>0$, a contradiction. Compressing $E$ itself
gives the upper bound, and the lower bound follows from~\eqref{eq-basic}.
\end{proof}

\begin{remark}\label{rem-ambient}
For any bounded set $E\subset\R^n$, the value $\Holdim E$ does not exceed the infimum
in~(i), which in turn does not exceed the infimum in~(ii). The three 
quantities coincide when $\Holdim E=\lfloor\dimH E\rfloor$. We do not know whether they coincide
in general. 
\end{remark}

\section{Finite upper box dimension}\label{sec-box}

We use the following embedding theorem of Foias and Olson~\cite[Theorem 1.2]{FoiasOlson}. See Robinson's monograph~\cite{Robinson} for more on such embeddings, and \cite[Corollary 2.6]{MargarisRobinson} for an explicit range of the exponent~$A$.

\begin{theorem}\label{thm-embedding}
Let $(X,d)$ be a metric space and let $N$ be an integer with $N>2\dimB X$. There exist a map $J\colon X\to\R^N$ and 
constants $c>0$ and $A\ge1$ such that 
\begin{equation}\label{eq-embedding}
c\,d(x,y)^A\le\|J(x)-J(y)\|_\infty\le d(x,y)\qquad \forall x,y\in X.
\end{equation}
\end{theorem}

The following elementary lemma is included for the sake of completeness.

\begin{lemma}\label{lem-box}
Finiteness of the upper box dimension is preserved by admissible changes
of metric.
\end{lemma}

\begin{proof}
If $\rho\le C\,d^\beta$, then every $d$-ball of radius $r$ is contained in
the $\rho$-ball of radius $Cr^\beta$ with the same center. Hence
$\dimB(X,\rho)\le\beta^{-1}\dimB(X,d)$, and the roles of $d$ and $\rho$ can
be interchanged.
\end{proof}

\begin{proposition}\label{prop-floor}
If $\dimB X<\infty$, then $\Holdim X\le\lfloor\dimH X\rfloor$.
\end{proposition} 

\begin{proof}
Apply Theorem~\ref{thm-embedding} to $X$ with an integer $N>2\dimB X$, and let $E=J(X)$ and $k=\lfloor\dimH X\rfloor$.
Since $J$ is Lipschitz, $\dimH E\le\dimH X<k+1$. For $\eps>0$ let $F$ be the map given 
by Theorem~\ref{thm-euclidean}. By~\eqref{eq-embedding} and~\eqref{eq-ambient-power}, the metric
$\rho(x,y)=\|F(J(x))-F(J(y))\|_\infty$ is admissible on $X$, and $\dimH(X,\rho)=\dimH F(E)\le k+\eps$.
\end{proof}

\begin{proof}[Proof of Theorem~\ref{thm-integer}]
The upper bound in~\eqref{eq-integer} is Proposition~\ref{prop-floor}, and
the lower bound is~\eqref{eq-basic}. Suppose that $k<\Holdim X<k+1$ for
some integer $k\ge0$, and choose an admissible metric $\rho$ on $X$ with
$\dimH(X,\rho)<k+1$. Since a composition of bi-H\"older maps is
bi-H\"older, $\Holdim(X,\rho)=\Holdim(X,d)$. By Lemma~\ref{lem-box},
$\dimB(X,\rho)<\infty$, and Proposition~\ref{prop-floor} applied to
$(X,\rho)$ gives
\[
\Holdim X=\Holdim(X,\rho)\le\lfloor\dimH(X,\rho)\rfloor\le k,
\]
a contradiction.

For the last assertion, let $\rho$ be an admissible metric on $X$. By
Lemma~\ref{lem-box}, $\dimB(X,\rho)<\infty$, and Theorem~\ref{thm-embedding}
applied to $(X,\rho)$ gives a Euclidean embedding of $X$ that is
bi-H\"older with respect to $d$ and whose image has Hausdorff dimension at
most $\dimH(X,\rho)$.
\end{proof}

\section{AI use disclosure}

The author used ChatGPT (OpenAI) and Claude (Anthropic) to assist with mathematical exploration and manuscript preparation. The author reviewed and verified the mathematical arguments, revised the manuscript, and takes full responsibility for its content.

\bibliographystyle{plain}

\bibliography{references-2.bib}

\begin{thebibliography}{10}

\bibitem{BMT}
Zolt\'an~M. Balogh, Roberto Monti, and Jeremy~T. Tyson.
\newblock Frequency of {S}obolev and quasiconformal dimension distortion.
\newblock {\em J. Math. Pures Appl. (9)}, 99(2):125--149, 2013.

\bibitem{BHW}
Christopher~J. Bishop, Hrant Hakobyan, and Marshall Williams.
\newblock Quasisymmetric dimension distortion of {A}hlfors regular subsets of a metric space.
\newblock {\em Geom. Funct. Anal.}, 26(2):379--421, 2016.

\bibitem{ColvinThesis}
Samuel Colvin.
\newblock {\em H{\"o}lder Dimension with Applications to Hyperbolic Groups}.
\newblock PhD thesis, University of Bristol, 2020.

\bibitem{Colvin}
Samuel Colvin.
\newblock Minimising {H}ausdorff dimension under {H}\"older equivalence.
\newblock {\em J. Lond. Math. Soc. (2)}, 103(3):781--816, 2021.

\bibitem{FoiasOlson}
Ciprian Foias and Eric~J. Olson.
\newblock Finite fractal dimension and {H}\"older-{L}ipschitz parametrization.
\newblock {\em Indiana Univ. Math. J.}, 45(3):603--616, 1996.

\bibitem{Heinonen}
Juha Heinonen.
\newblock {\em Lectures on analysis on metric spaces}.
\newblock Universitext. Springer-Verlag, New York, 2001.

\bibitem{HurewiczWallman}
Witold Hurewicz and Henry Wallman.
\newblock {\em Dimension {T}heory}, volume~4 of {\em Princeton Mathematical Series}.
\newblock Princeton University Press, Princeton, NJ, 1941.

\bibitem{Kovalev}
Leonid~V. Kovalev.
\newblock Conformal dimension does not assume values between zero and one.
\newblock {\em Duke Math. J.}, 134(1):1--13, 2006.

\bibitem{MackayTyson}
John~M. Mackay and Jeremy~T. Tyson.
\newblock {\em Conformal dimension}, volume~54 of {\em University Lecture Series}.
\newblock American Mathematical Society, Providence, RI, 2010.
\newblock Theory and application.

\bibitem{MargarisRobinson}
Alexandros Margaris and James~C. Robinson.
\newblock Embedding properties of sets with finite box-counting dimension.
\newblock {\em Nonlinearity}, 32(10):3523--3547, 2019.

\bibitem{Pansu}
Pierre Pansu.
\newblock Dimension conforme et sph\`ere \`a{} l'infini des vari\'et\'es \`a{} courbure n\'egative.
\newblock {\em Ann. Acad. Sci. Fenn. Ser. A I Math.}, 14(2):177--212, 1989.

\bibitem{Robinson}
James~C. Robinson.
\newblock {\em Dimensions, embeddings, and attractors}, volume 186 of {\em Cambridge Tracts in Mathematics}.
\newblock Cambridge University Press, Cambridge, 2011.

\bibitem{Romney}
Matthew Romney.
\newblock Singular quasisymmetric mappings in dimensions two and greater.
\newblock {\em Adv. Math.}, 351:479--494, 2019.

\bibitem{Salem}
Rapha\"{e}l Salem.
\newblock On some singular monotonic functions which are strictly increasing.
\newblock {\em Trans. Amer. Math. Soc.}, 53:427--439, 1943.

\bibitem{Tukia}
Pekka Tukia.
\newblock Hausdorff dimension and quasisymmetric mappings.
\newblock {\em Math. Scand.}, 65(1):152--160, 1989.

\end{thebibliography}

\end{document}